\input amstex
\documentstyle{amsppt}
\magnification = \magstep 1
\NoBlackBoxes
\vsize=9truein
\hsize=6.5 true in

\define\({ \left (}
\define\){\right )}

\define\brq{^{[q]}}
\define\brt{^{[t]}}
\define\q{^{1/q}}
\redefine\vec#1 #2{#1_1,\ldots, #1_{#2}}

\define\bx{\bold x}
\define\bfy{\bold y}
\define\bz{\bold z}

\define\m{\bold m}
\define\n{\bold n}
\define\um{\underline{\m}}
\define\un{\underline{\n}}

\define\inc{\subseteq}
\define\1{^{-1}}
\define\8{{\infty}}
\define\Ro{R^{\circ}}
\redefine\Re{R^{(e)}}
\define\defn#1{{\sl #1}}

\define\bd{\bold}

\define\Ga{\Gamma}
\define\hatR{{\widehat R}}
\define\hatS{{\widehat S}}
\define\Shat{{\widehat S}}
\define\Rhat{{\widehat R}}
\define\ttau{{\widetilde \tau}}

\define\al{\alpha}
\define\eps{\epsilon}
\define\la{\lambda}

\def\Ann {\operatorname{Ann}\nolimits}
\def\Hom {\operatorname{Hom}\nolimits}

\def\im {\operatorname{Im}\nolimits}
\def\hgt {\operatorname{ht}\nolimits}

\def\Spec {\operatorname{Spec}\nolimits}

\def\ord {\operatorname{ord}\nolimits}
\def\exp {\operatorname{exp}\nolimits}
\topmatter

\title \nofrills{Test ideals and base change problems in tight closure theory} 
\endtitle
\author Ian M. Aberbach and
Florian Enescu
\endauthor
\leftheadtext{Aberbach, Enescu, Test ideals and base change problems}
\rightheadtext{Aberbach, Enescu, Test ideals and base change problems}
\address
Department of Mathematics, University of Missouri,
Columbia, Missouri 65211
\endaddress
\email aberbach\@math.missouri.edu \endemail
\address
Department of Mathematics, University of Michigan,
Ann Arbor, Michigan  48109 and
Institute of Mathematics of the Romanian Academy,
Bucharest, Romania
\endaddress
\email
fenescu\@umich.edu
\endemail
\thanks
The first author was partially supported by the NSF and by the University
of Missouri Research Board. The second author thanks the University of Michigan for support through the Rackham Predoctoral Fellowship.
\endthanks
\subjclass Primary 13A35, Secondary 13B40
\endsubjclass
\abstract Test ideals are an important concept in tight closure theory and their behavior via flat base change can be very difficult to understand. Our paper presents results regarding this behavior under flat maps with reasonably nice 
(but far from smooth) fibers. This involves analyzing, in depth, a special type of ideal of test elements, called the CS test ideal. Besides providing new results, the paper also contains extensions of a theorem by G.~Lyubeznik and K.~E.~Smith on the completely stable test ideal and of theorems by F. Enescu and, independently, M. Hashimoto on the behavior of $F$-rationality under flat base change.

\endabstract

\endtopmatter

\document
\head 1. Introduction and terminology \endhead

Let  $R$  be a    Noetherian  commutative  ring   of positive  prime
characteristic $p$.   Over the last  decade  tight closure  theory has played  
a tremendous  role  in understanding    the structure of  $R$, \cite{9--12},
\cite{16}.  One compelling problem is how tight closure behaves under flat base
change.   This  problem has multiple facets and   it is intimately connected 
to  the problem of localization of tight closure, as shown by M.  Hochster and
C. Huneke  in \cite{11}.  In fact, they show that even the case when  the
fibers are ``good'' (as in the case   of   smooth base change)    is very
complicated and  not entirely understood.  Several  important questions remain 
open even in the case of the completion morphism   of a local excellent ring.
With this  in mind, the  goal of understanding how tight closure behaves under
flat base  change, while weakening the smoothness conditions on the fibers, is quite challenging. 

The key idea has become to allow a wider class of fibers.
Recently  techniques
have been developed to better
understand how $F$-rationality, weak $F$-regularity, and
strong $F$-regularity  behave under   certain flat   local homomorphisms $(R,
\m) \to (S,  \n)$  with non-smooth fibers (\cite{3}, \cite{7}, \cite{8}). 
It is now  apparent that
significant progress can be obtained in  studying the tight closure of
$\m$-primary ideals for flat local homomorphisms (as the results in \cite{3}
indicated). At the same  time, Lyubeznik and  Smith  initiated a new  approach
to handle the localization and completion problem for certain test ideals
\cite{16}. Test elements are a crucial tool  in tight closure theory  and the
question of  how they behave  under  localization  is  one of the  central   
problems of the theory. 
In fact, the knowledge that $1 \in R$ is a test element is {\it not} (yet)
known to imply that $1 \in R_P$ is a test element!
Lyubeznik and Smith suggested that one should
understand first the case of  a particular ideal contained in the  test ideal
(and  conjecturally equal to  it in the case of  excellent reduced rings).  In
the local reduced case, this ideal, $\ttau(R)$, is the annihilator of the tight
closure of $0$ in the injective hull of the residue field. They  show that 
$\ttau(R)$ (called here {\it the $CS$  test ideal}) behaves as expected with
respect   to localization and  completion  for   rings that are images  of
excellent regular local rings. 
Although able to prove the preceding statement,
 Lyubeznik and Smith  expected this to be true for {\it all} excellent
local rings (see the first paragraph after Theorem 7.1 in \cite{16}).

In our paper we investigate the behavior of tight closure of ideals primary to
the maximal ideal (in local excellent reduced rings)  under flat base change.
Our work relies on the CS test ideal and its properties. The conditions on
the fibers of the homomorphism are greatly  weakened and our results apply  in
many situations other than that of ideals primary to the maximal ideal. We
explore this in section 4.  However, as we just said, we must develop several
theorems about test elements which, in fact, are interesting in their own
right.  These appear in section 3.  Along the way we show that the CS
 test ideal respects both completion and localization in excellent 
reduced local rings, as Lyubeznik and Smith expected. 
 Extending their result to the class of excellent rings is more than a technical achievement. On one hand, J. Nishimura has given examples of excellent rings that are not images of regular rings (see \cite{18}). On the other hand, 
while much of tight closure theory has been shown to apply well to the class of excellent rings, recent examples of 
Loepp and Rotthaus show that if one goes beyond this class, then tight closure does not even commute with completion, \cite{15}. 
Our theorem provides more evidence in support of the belief that tight 
closure behaves well in excellent rings. To prove
the theorem we make use of Lyubeznik's and Smith's result that the CS test ideal localizes in complete
local reduced rings. The reduction is done by the very means of flat base change, yet in a non-canonical manner. It underlines the theoretical importance of results on the behavior of tight closure under such base change and
highlights how understanding even smooth base change depends on understanding
flat base change where the fibers are {\it not} regular.

This result appears below as:
\proclaim{Theorem 3.6} Let $(S,\un)$ be a 
 semi-local reduced excellent ring.
Then
\roster
\item $\ttau(S)$ localizes properly at any prime $Q \in \Spec(S)$,
i.e., $\ttau(S)_Q = \ttau(S_Q)$.
\item $\ttau(\hatS) = \ttau(S) \hatS$. 
\item 
 Let $(R,\m)$ be a  semi-local reduced excellent ring
 such that $(R,\um)\to (S,\un)$ is a flat semi-local
map 
with Gorenstein $F$-injective closed fibers and
such that if $\n$ is a maximal ideal of $S$ lying over $\m$, then
$R/\m \to S/\n$ is separable.  
 Assume that $R$ and $S$ have a common CS 
 test element.
 Then $\ttau(S) = \ttau(R)S$.
\endroster
\endproclaim

From this we obtain that the CS test ideal extends properly in a number
of good cases:
\proclaim{Corollary 3.8} Let $(R,\m) \to (S,\n)$ be a flat
map of reduced  local rings such that $S/\m S$ is
Gorenstein and $0^*_{E_S}$ is extended.   Then
$\ttau(S) = \ttau(R)S$.
\endproclaim

It is clear that test elements are important in problems of base change. Our investigation also shows that it 
suffices to have available a particular type of test elements, namely the CS
test  elements. Generally, these elements  are as abundant as the test elements
and one now has the advantage that the theory of the CS test ideal  is well
understood (as our results just quoted above, together with those of Lyubeznik and Smith, show).

In section 4 we turn to the question of extending tight closure
by flat base change. One of the main cases where our results apply 
is that of  ideals primary to the maximal ideal in local excellent rings.
We extend many of the known results.  To be precise, Hochster and Huneke (\cite{11}) have results on the behavior of tight closure of ideals under smooth (or regular, using a more often used terminology) base change; we extend their results by weakening the conditions on the fibers. Aberbach (\cite{3}) has results on the behavior of weak $F$-regularity obtained by considering the behavior of irreducible ideals; in contrast, we obtain results for all ideals  primary to the maximal ideal. 
These results are considerably more difficult, since non-irreducible ideals
have larger dimension socles, which turn out to require new ideas for
carrying out tight closure arguments.
Enescu (\cite{7}) and, independently,  Hashimoto (\cite{8}) have results on the behavior of $F$-rationality under some finiteness conditions (their theorems in turn extend V\'elez's results on $F$-rationality, \cite{20}); ours show that these conditions are not necessary. A feature of the study is that some of the tools on which we rely heavily are as basic as Nakayama's Lemma.

We would like to illustrate some of our results here. For example, Theorem 4.1 shows that tight closure
extends well when the closed fiber is Gorenstein (and either $F$-injective
or $F$-rational). However, when the closed fiber is assumed to be only
CM, the arguments have been more difficult.  We have been able to
show:
\proclaim{Theorem 4.2}  Let $(R,\m,K) \to (S,\n,L)$ be a flat map
of excellent rings with CM closed fiber.  Let $I \inc R$
be $\m$-primary and tightly closed.  Let $\bz = \vec z d$ be 
elements whose image in $S/\m S$ is an s.o.p.  Suppose that either
\roster
\item the element $c \in \Ro$ is a common test element for $R$ 
and $S$, and that $S/\m S$ is geometrically 
 $F$-injective over $K$, or
\item $S/\m S$ is geometrically $F$-rational over $K$.
\endroster
Then $(I,\bz)S$ is tightly closed.
\endproclaim

\comment
Let us describe now a byproduct of our work that addresses the notion of test elements.
When $R \to S$ is a ring homomorphism, there are a number of results
that involve the need for a {\it common} test element.  This is a
somewhat restrictive requirement, and even when there is such a
common test element, its existence may not be known a priori.
To this end we have proved several results giving CS test elements
for $S$ having a convenient form (relative to $R$).  It can be argued that these results are mainly of technical significance. However, much of our results 
in section 4 makes heavily use of them. We believe that they will prove to be useful in contexts diferent from ours and would like to highlight one of them here. As its proof shows (see section 5), the theorem is strongly related to the behavior 
of CS test ideal under flat maps with good fibers. We succeed keeping a part of the desired behavior of test elements that is significant in applications intact, while replacing the hypotheses of Theorem 3.6 (3) with conditions on a single fiber.

\proclaim{Theorem 5.6} Let $(R,\m) \to (S,\n)$ be a map of
excellent reduced rings.  Assume that $\Rhat$ is
a domain.  Let $Q \in \Spec(S)$ be such that $\m S \inc Q$,
$R \to S_Q$ is flat, and $\m S_Q$ is reduced.  Assume that either
\roster
\item
$\widehat {S_Q}$ is a domain and $\widehat {S_Q}/ \m \widehat {S_Q}$
is Gorenstein and $F$-rational, or
\item
 $S_Q/\m S_Q$ is regular.
\endroster
Then there is an element $d \in S -Q$ such that $d \ttau(R) \inc \ttau(S)$.
\endproclaim

\endcomment

When $R \to S$ is a ring homomorphism, there is a large number of results in the literature
that involve the need for a {\it common} test element (see for example our Theorem 3.3).  This is a somewhat restrictive requirement, and even when there is such a
common test element, its existence may not be known a priori. Section 5 deals with this issue. We use the notion of strong test ideals to show how this assumption can be removed in certain conditions (that are more general than the ones in~\cite{11}, see for example Discussion 7.11) and still have that tight closure behaves well: 

\proclaim{Theorem 5.3} Let $(R,\m, K) \to (S,\n, L)$ be a flat map of
 excellent domains both of whose 
completions are domains.  Assume that $S/\m S$ is Gorenstein and $F$-rational
and $K \to L$ is separable.  Then $0^*_{E_S}$ is extended from $0^*_{E_R}$,
and $\ttau(S) = \ttau(R)S$.

Also, $0^{*fg}_{E_S}$ is extended from $0^{*fg}_{E_R}$ and hence $R$ and
$S$ have a common test element.
\endproclaim

In fact, we prove results that show the existence of CS test elements for $S$ of a convenient form relative to $R$ (Theorems 3.9 and 5.6). It can be argued that they are mainly of technical significance. However, many of the theorems in section 4 make use of such ideas. We believe that they will prove to be useful in a number of contexts different from ours.

In the last section of the paper we give results concerning Hilbert-Kunz
multiplicities and test exponents.

We now define tight closure and related concepts (see \cite{9} and
\cite{11}).  We use $q$ to denote a power of $p$; so $q = p^e$ for $e \ge 0$. 
For $I \inc R$ set $I\brq = (i^q: i\in I)$.  Let $\Ro$ be the complement in $R$
of the minimal primes of $R$.  We say that $x \in I^*$, the tight closure of
$I$, if there exists $c \in \Ro$ such that for all $q \gg 0$, $c x^q \in I\brq$.
We say that $x$ is in the \defn{Frobenius closure} of $I$, $I^F$, if there
exists a $q$ such that $x^q \in I\brq$, and say that $I$ is \defn{Frobenius
closed} if $I = I^F$. When $R$ is reduced, then $R\q$ denotes the ring of $q$th
roots of elements of $R$.  Note that in this case, $cx^q \in I\brq$ if and only
if  $c\q x \in IR\q$.  When $R\q$ is module-finite over $R$, we call $R$
\defn{$F$-finite}.  We call $R$ \defn{weakly $F$-regular} if every ideal of $R$
is tightly closed.  A weakly $F$-regular ring is always normal, and under mild
hypotheses is Cohen-Macaulay.  $R$ is \defn{$F$-regular} if every localization
of $R$ is weakly $F$-regular.  An $F$-finite reduced ring $R$ is called
\defn{strongly $F$-regular} if for all $c \in \Ro$ there exists $q$ such that
$Rc\q \inc R\q$ splits over $R$.  Strongly $F$-regular implies $F$-regular.  It
is open in most cases whether or not these notions are equivalent.

We call an ideal $I = (\vec x n)$ a \defn{parameter ideal} if $\hgt I \ge n$.
$R$ is \defn{$F$-rational} if every parameter ideal is tightly closed.
$F$-rational and Gorenstein rings are (strongly) $F$-regular.  We call
the local ring $(R,\m)$ \defn{$F$-injective} if the Frobenius endomorphism
$R \to R$ sending $r \mapsto r^p$ induces an injection on all the
local cohomology modules of $R$.  When $R$ is CM, then $R$ is $F$-injective
if and only if some (equivalently, every) ideal generated by a 
system of parameters is Frobenius closed.  
\comment For rings which are not
necessarily CM we will use the term \defn{parameter $F$-pure} to denote
the condition that every s.o.p.~is Frobenius closed (this may be
weaker than $F$-injective).  
\endcomment

There is a tight closure
operation for a submodule $N \inc M$.  Let $\Re$ be the $R$-bimodule whose
underlying group structure is the same as for $R$, considered as a left
$R$-module in the usual way, but as a right $R$-module via the $e$th iteration
of the Frobenius map sending $r \mapsto r^{p^e}$.  Then $F^e_R(M) = \Re
\otimes_R M$ considered as a left $R$-module.  For $m \in M$, we
let $m^q = 1 \otimes m \in F^e_R(M)$.  If $N \inc M$, then we
let $N\brq_M$ be the image of $F^e_R(N)$ in $F^e_R(M)$.  An element $x \in M$
is in the tight closure of $N$ in $M$, denoted $N^*_M$ if there exists 
$c \in \Ro$ such that for all $q \gg 0$, $cx^q \in N\brq_M$.  
The \defn{finitistic tight closure of $N$ in $M$}, denoted $N_M^{*fg}$ is
the union of all $(N \cap M')^*_{M'}$ where the union is taken over
all finitely generated submodules $M' \inc M$.  
 
We will often be considering modules of the form
$M = \varinjlim_t R/I_t$ for a sequence of ideals
$\{I_t\}$.  Thus it is useful to discuss tight 
closure in this case.  Let $u \in M$ be an element which is given by
$\{u_t\}$ where in the direct limit system $u_t \mapsto u_{t+1}$.
Then $u \in 0^*_M$ if there exists $c \in \Ro$
and a sequence $t_q$ such that for all $q \gg 0$,
$c u_{t_q}^q \in I_{t_q}\brq$,  and $u$ is in  $0^{*fg}_M$ if
there exists $c \in \Ro$ and $t >0$ such that $c u_t^q
\in I_t\brq$ for all $q$.  Clearly $0^{*fg}_M \inc 0^*_M$.

The element $c \in \Ro$ is a \defn{test element} if whenever $x \in I^*$ (for
any $I$), then $cx^q \in I\brq$ for all $q$.  The element $c \in \Ro$ is a
\defn{completely stable test element} if $c$ is a test element in every
completion of every localization of $R$.  When $R$ is an excellent reduced
local ring, then $R$ has an abundance of test elements.  When $R$ is excellent
and reduced, which is the only case we will consider here, then any (completely
stable) test element for ideals is also a (completely stable)
test element for tight closure in
finitely generated modules.  Of course, we will be very interested in
tight closure in non-finitely generated modules such as the injective
hull of the residue field.  We deal with this situation in section 2.

\proclaim{Theorem 1.1} {\rm (\cite{11}, Theorem 6.2)}  Let $(R,m)$ be a reduced
excellent local ring.  If $c \in \Ro$ and $R_c$ is Gorenstein and 
$F$-regular, then $c$ has a power which is a completely stable
test element.
\endproclaim

One way to understand the singularities of $R$ is to understand
the \defn{test ideal} (respectively, the \defn{completely stable test
ideal}).  The test ideal of $R$ is $\tau(R) = \{c \in R: \text{for all 
$I$ and all $q$, $x \in I^*$ if and only if $cx^q\in I\brq$}\}$. The completely stable
test ideal  is defined analogously.  Suppose that $(R,\m)$ is
reduced and approximately Gorenstein, i.e., there exists a sequence
$\{I_t\}$
of $\m$-primary irreducible ideals cofinal with the powers of $\m$.
Then $\tau(R) = \cap_t (I_t:I^*_t)$.  If $E_R = E_R(R/\m)$ is the injective
hull of the residue field, then we also have $\tau(R) = \Ann_R 0^{*fg}_{E_R}$.

We would like to thank the referee for useful editorial suggestions and Graham Leuschke for providing us a copy of the preprint~\cite{18}.

\head 2. CS Test elements 
\endhead

 As mentioned earlier,
Lyubeznik and Smith have defined the ideal $\ttau(R) = \Ann_R 0^*_{E_R}$,
which is contained in the completely stable test ideal (and is probably
equal to it when $R$ is a reduced excellent local ring).  

\definition{Definition}
We will call
$\ttau(R)$ the \defn{CS test ideal} of $R$. The elements of $\ttau(R)$
that are in $\Ro$ will be called \defn{CS test elements}.
\enddefinition

Remember that $\tau(R) = \Ann_R 0^{*fg}_{E_R}$.
It is conjectured, but open in most cases, that in fact $0^*_{E_R} =
0^{*fg}_{E_R}$.  Clearly, $\ttau(R) \inc \tau(R)$. 

When $R$ is not local, we will use the convention in
\cite{16} that $E_R = \oplus_\m E_{R_\m}$ where the sum is over all
maximal ideals of $R$.  It is natural to ask if given a multiplicative
system $U \inc R$, we have $\tau(U\1 R) = U\1 \tau(R)$ (or if
$\ttau(U\1 R) = U\1 \ttau(R)$).   The latter equality is obtained
for $R$, a reduced ring that is a homomorphic image of an excellent regular ring as in \cite{16}.
We show that in fact this equality holds for all excellent reduced local
rings (see Theorem 3.6).   Understanding localization of $\tau(R)$ is actually
a more difficult task since knowing $\tau(R)$ localizes will imply that
weak $F$-regularity is equivalent to $F$-regularity.

As mentioned above, understanding the (finitistic) tight closure of $0$
in $E =E_R(R/\m)$ for the local ring $(R,\m)$ is very important.  When
$(R,\m)$ is excellent and reduced, then $E = \varinjlim_t R/I_t$ where
the $I_t$ are irreducible $\m$-primary ideals cofinal with the powers
of $\m$.  Thus the (finitistic) tight closure of $0$ in $E$ may be
understood as described above.

Proposition 2.1 below shows that CS test elements are precisely those
elements that test tight closure in Artinian modules.
  Without loss of generality, it is
enough to show that all such elements need work only in tight closure tests for
$0^*_M$.

\proclaim{Proposition 2.1}  Let $(R,\m,k)$ be an excellent reduced ring, and
M an Artinian $R$-module.  Then
\roster
\item $F^e(M)$ is Artinian for all $e$, and 
\item if $c \in \ttau(R) \cap \Ro$, then $c$ works in all tight closure 
tests for $0^*_M$.
\endroster
\endproclaim

\demo{Proof}  (1)  We may assume that $R$ is complete.  We can then choose
$\Ga$ such that $S = R^\Ga$ is still reduced.  Note that an $R$ module
$N$ is Artinian over $R$ if and only if $N \otimes_R S$ is Artinian over $S$.
Also, $F^e_R(M) \otimes_R S = F^e_S(M \otimes_R S)$; so without loss of
generality we may assume that $R$ is $F$-finite.  But in general, if
$R \to S$ is module finite, then $N$ is Artinian over $R$ if and only if
$N \otimes_R S$ is Artinian over $S$.  In this case, we are taking 
$S = R^{(e)}$.

(2)  Suppose that $M$ is Artinian and $u \in N^*_M$.  Since $M/N$ is Artinian,
we may assume that $N =0$.  Then $u^q \in 0^*_{F^e(M)}$ for all $e$.
By part (1), $F^e(M)$ is Artinian, hence $F^e(M) \inc E^{n_e}$ (where $E = E_R(k))$. (To see this it is enough to tensor with $\hatR$.) The image
of $u^q$ is in the tight closure of $0$ in $E^{n_e}$.  Thus
$cu^q = 0$. \qed
\enddemo

\remark{Observation}  The argument in Proposition 2.1 shows that, in fact,
elements of $\ttau(R)$ test tight closure in {\it any} module $M$ supported
only on $\{\m\}$.
\endremark

The test elements shown to exist in Theorem 1.1 are actually CS test elements.
We refer the reader to section 6 of \cite{11} for a full development
of $\Ga$-extensions.  The main point is that if $R$ is a complete local
ring, then one may extend the coefficient field to obtain a 
map $R \to R^\Ga$ which is faithfully flat and purely inseparable
and such that $R^\Ga$ is $F$-finite.  Moreover, if $R$ is reduced, then
$\Ga$ may be chosen so that $R^\Ga$ is reduced.

\proclaim{Lemma 2.2}  Let $(R,\m,k)$ be an excellent reduced ring.
If $R_c$ is Gorenstein and weakly $F$-regular then $c$ has a
power which is a CS test element.
\endproclaim

\demo{Proof} We already know that $c$ has a power (which we will rename as $c$)
which is a test element by Theorem 1.1.  In fact, the proof, which involves 
looking at $\hatR^\Ga$ for sufficiently small $\Ga$, along with Proposition
2.1, shows the CS property.  Because $\Rhat^\Ga$ is faithfully flat over $R$,
it is enough to consider the case when $R$ is $F$-finite. If $E_R(k) =
\varinjlim_t R/I_t$ and $u = \{u_t\} \in 0^*_E$, then there exists $d \in \Ro$
such that for all $q$ there exists $t$ such that $d u_t^q \in I_t\brq$. For 
$R$ there is a $q_0$ and a linear map  of $R^{1/q_0} \to R$ sending  $d^{1/q_0}
\to c$. Since  $du_t^{qq_0} \in I_t^{[qq_0]}$  we may take $q_0^{th}$ roots and
apply the splitting to see that 
$cu_t^q \in I_t^{[q]}$. Thus $c \in \Ann_R 0^*_E$.
\qed \enddemo

\proclaim{Theorem 2.3} {\rm (see \cite{11}, Discussion 7.11)}
Let $(R,\m,k) \to (S,\n,L)$ be a flat map of complete local rings
with regular closed fiber.  Assume also that $k \to L$
is separable.  Then $R$ and $S$ have a common CS test
element.
\endproclaim

\demo{Proof} In  light of  Lemma 2.2, it is enough to find an element $c \in R$
such that $S_{c}$ is regular. We will follow  7.11 in \cite{11}. As it is done there,
take $K$ to be a coefficient field for $R$, $\bx = x_{1},...,x_{n}$ a system of
parameters in $R$ such that $R$ is module finite over $A = K[[\bx]]$. Let $\bfy
= y_{1}, ..., y_{d}$ be a regular system of parameters for $S/\m S$ and enlarge
$k$ to a coefficient field $L$ of $S$  (we can do this because $L$ is
separable over $k$). It can be shown, as in the reference
quoted above, that $S = B \otimes_{A} R$ where $B = L[[\bx,\bfy]]$. Now it is
clear that if $c \in R$ is such that $R_{c}$ is regular, then $S_{c}$ is
regular (the base change from $R$ to $S$ has regular fibers). Now, 
Lemma 2.2 applies and we have a common CS test element for $R$ and $S$. \qed
\enddemo

We end this section by stating a result that relates $E_S$ to $E_R$
for a flat local homomorphism $(R,\m) \to (S, \n)$ with Gorenstein closed fiber.
Our results rely heavily on this structure theorem.

\proclaim{Theorem 2.4} {\rm (\cite{11})} Suppose that $(R, \m) \to (S, \n)$ is
a flat local homomorphism such that $S/\m S$ is Gorenstein. If $\bz$ is a
regular sequence for the closed fiber, then $E_S = \varinjlim_t E_R \otimes_R
S/(\bz)\brt$. \endproclaim

\comment 
The following theorems have essentially the same proofs as
those given in \cite{9}.

\proclaim{Theorem 2.2} {\rm \cite{9, Theorem 6.9}}  Let $R$
be  torsion-free, module-finite, and generically smooth
over a regular local domain $(A,\n)$.  Let $c \in A^\circ$
be such that $c R^\infty \inc A^\infty[R]$.  Then
the following conditions for an element $x = \{x_t\} \in 
M = \varinjlim_t R/I_t$ are equivalent:
\roster
\item  $x \in 0^*_M$.
\item  $x \in 0^*_{M'}$ over $\hatR$ where $M' = \varinjlim_t
\hatR/I_t \hatR$.
\item  There exists a sequence of elements $\{\eps_n\}$ in
$((\hatR)^\8)^\circ$ such that $\ord(\bold N(\eps_n)) \to 0$
as $n \to \8$ and for infinitely many $n$ there
exists $t$ such that $\eps_n x_t \in I_t (\hatR)^\8$.
\item $cx^q =0$ in $F^e(M)$ for all $q$.
\endroster
\endproclaim

\proclaim{Theorem 2.3} {\rm \cite{9, Theorem 6.13}} Let $(R,\m)$
be module-finite, torsion-free, and generically smooth
over a regular local domain $(A,\n)$.  Then every element $d \in A^\circ$
such that $R_d$ is smooth over $A_d$ has a power $c$ which is
a DL completely stable test element in $B \otimes_A R$ for every $A$-flat
regular domain such that $B \otimes_A R$ is local.  A sufficient
condition for $c$ to have this property is that 
$cR^\8 \inc A^\8[R]$.
\endproclaim

\proclaim{Theorem 2.4} {\rm \cite{9, Corollary 6.17}}  Let $(R,\m)$
be a reduced local ring that is module-finite and torsion-free
over an excellent local regular ring $(A,\n)$.  If $q'$
is sufficiently large so that $R[A^{1/q'}]$ is generically smooth
over $A^{1/q'}$ and $c \in A^\circ$ is such that $c^{1/q'}R^\8
\inc A^\8[R]$ (such elements always exist), then $c$ is a
DL completely stable $q'$-weak test element in $B\otimes_A R$ for
every local excellent regular $A$-flat domain $B$ such that
$B\otimes_A R$ is local and reduced.
\endproclaim

\proclaim{Theorem 2.5} {\rm \cite{9, Corollary 6.19}}
Let $(R,\m)$ be a reduced, equidimensional, complete local
ring.  Then $R$ has a DL completely stable $q'$-weak test
element $c$.
\endproclaim

We may also get an extended version of a result from \cite{HH?}
on common test elements:

\proclaim{Theorem 2.6} {\rm \cite{9, Corollary 6.?}}
Let $(R,\m) \to (S,\n)$ be a flat map of complete local rings
with regular closed fiber.  Assume also that $R/\m \to S/\n$
is separable.  Then $R$ and $S$ have a common DL weak test
element.
\endproclaim

\endcomment

\head \S 3.  Completely stable test ideals under flat base change \endhead

The main theorem in this section, Theorem 3.6,
 shows that if $(R,\m) \to (S,\n)$ 
is a flat map of excellent reduced rings
with sufficiently nice closed fiber, then the CS test
ideal of $S$ is extended from $R$. It also shows that  for any
reduced excellent (semi-)local ring, the completely stable test
ideal localizes properly.  

  We also wish to highlight Theorem 3.9.  There are a number of results in the
literature concerning a (flat) map $R \to S$ where one assumption is that $R$
and $S$ have a common test element.  This assumption can be somewhat
restrictive.  Theorem 3.9 shows that $S$ can have useful test elements in
relation to $R$ even if there is no common test element.  We will make
extensive use of this result in section 4.

In the remainder of this paper we will make frequent use of the
following facts \cite{17}:
\roster
\item  If $R \to S$ is flat,   $N \inc M$ are finitely generated
$R$-modules, and $w \in M$, then $(N \otimes S):_S (w \otimes 1) = (N:_R w)
\otimes S$.
\item If $(R,\m) \to (S,\n)$ is a flat local map of local rings and
$\bz = \vec z d \in S$ is a regular sequence on $S/\m S$, then
$R \to S/(\bz)S$ is flat.
\endroster

\proclaim{Lemma 3.1}  Let $(R,\m) \to (S, \n)$ be a flat map of reduced local
rings with $F$-injective CM closed fiber and a common CS test element $c$. 
Suppose $\bz = \vec z d$ is an s.o.p. in $S/\m S$ and $b$ maps to
a socle element in $S/(\m, \bz)S$.  If $M$
 is an Artinian module and $u$ is not in $0^*_M$,
then $bu$ is not in $0^*_{M \otimes S/(\bz)}$.
\endproclaim

\demo{Proof}   Let $M' = M \otimes S/(\bz)$.  If the lemma is not true,
then for all $q$,  $c(bu)^q =0$ in $F^e_S(M')$.
Thus $b^q \in 0 :_S cu^q = (0 :_R c u^q)S + (\bz)\brq S$,
by flatness.  Since $u \notin 0^*_M$, $(0 :_R c u^q) \inc \m$
for infinitely many $q$.  Thus $b^q \in (\bz)\brq (S/\m S)$, contradicting
the $F$-injectivity of $S/\m S$. \qed
\enddemo

\remark{Remark 3.2}
Let $(A,\m) \to (B, \n)$ be a flat map of local rings such that
$\dim A = 0$.  The map is then faithfully flat.  In this situation
$B$ is a free $A$-module and any basis of $B/\m B$ lifts to a free
basis of $B$ over $A$.  To see this, let $\{\overline b_\la\}_{\la
\in \Lambda}$ be a basis for $B/\m B$ and map a free $A$-module $G
= \oplus_{\la \in \Lambda} A$ to $B$.  Then we have $B = \im(G) + \m \im(G)$
and since $\m$ is nilpotent, Nakayama's Lemma gives $B = \im(G)$
\cite{17}.  If $K = \ker(G \to B)$, then we see that $K \otimes A/\m
= 0$ by flatness of $B$ over $A$, and so $K=0$ by another application
of Nakayama's Lemma.
\endremark

If $\bz = \vec z d$ are elements of a ring $S$ then we will use
$(\bz)\brt$ to denote the ideal $(z_1^t,...,z_d^t)S$ for $t \ge 1$.

\proclaim{Theorem 3.3}  Let $(R,\m,K) \to (S,\n,L)$ be a flat
map of excellent rings
with Gorenstein $F$-injective closed fiber.  Assume that $K \to L$ is
a separable field extension.  Let $M = \varinjlim R/I_t$ be an
Artinian module which can be written as a 
direct limit of $\m$-primary ideals.  Suppose that $0$ is tightly closed
in $M$,
 and let $\bz = \vec z d \in S$ be elements
in $S$ which form an s.o.p.~in $S/\m S$.  Suppose that
 $c$ is a common CS test element for $R$ and $S$.
Then $0$ is tightly closed in $M' = \varinjlim S/(I_t + (\bz)\brt)S$.

In particular, if $R$ and $S$ are
reduced, then $0^*_{E_S}$ is the image of $0^*_{E_R}$ in $E_S$.
\endproclaim

\demo{Proof}  
Let $b \in S$ map to the socle of $S/(\m + (\bz))S$.  
For any element $u = \{u_t\} \in M$, suppose that for all $q$ there exists $t$
such that $c (u_t(z_1 \cdots z_d)^{t-1}b)^q \in (I_t, z_1^t,\ldots, z_d^t)\brq$.
Then, since the fiber is CM, we have $c (u_t b)^q \in (I_t, \bz)\brq$.

Given any 
choice of $\vec u n$ generating the socle of $M$, the images of the set
$\{u_1 b,\ldots, u_n b\}$ generate the socle of $M'$.
If $0$ is not tightly closed in $M'$, then there is a socle element
$v = \sum_{i=1}^n s_i u_i b$ in $0^*_{M'}$, and the $s_i$ can
be considered to be elements of $L$.  Choose $\vec u n$ so that such 
a sum can be written with the fewest possible $u_i$'s.  Let $h$ be
this number and reorder so that $v = \sum_{i=1}^h s_i u_i b$.
Suppose that the elements $\vec s h$ (thought of as elements of $L$)
are linearly dependent over $K$.  Then without loss of generality,
$s_h = r_1 s_1 + \cdots + r_{h-1} s_{h-1}$ and hence 
$v =  \sum_{i=1}^{h-1} s_i(u_i + r_i u_h)b$,
contradicting our choice of $h$.  Also, $h \ne 1$, by Lemma 3.1.

  We claim that for infinitely many
 $q$ there exist $t$ (depending on $q$) such that the
set $\{(s_1b)^q,\ldots, (s_hb)^q\}$ is part of a free basis over
$R/I_t\brq$ for  $S/(I_t,\bz)\brq S$.  If not, then by  Remark 3.2,
for infinitely many $q$ there exist $r_{iq}$
not all in $\m$ such that $\sum_{i=1}^h r_{iq} s_i^q b^q \in \m S + (\bz)\brq 
S$.  Since $L$ is separable over $K$, each $\sum_{i=1}^h r_{iq} s_i^q$
is a unit (or else, modulo $\n S$, we have a relation on
$\vec s h$ over $K^{1/q}$)
 and thus $b^q \in \m S + (\bz)\brq S$, contradicting
the $F$-injectivity of $S/\m S$.  Now, if $cv^q = \sum_{i=1}^h
(cu_i^q)(s_ib)^q =0$ in $F^e(M')$ for infinitely many $q$, we must have
$cu_i^q =0$ in $F^e(M)$ for infinitely many $q$
(by the assertion about free bases above).  
Thus $u_i \in 0^*_M$, a contradiction.

The last statement follows by taking $M = E_R/0^*_{E_R}$.  When
$R$ is excellent and reduced, it is approximately Gorenstein, and
$E_R$ may be written as a direct limit of the desired form.
By Theorem 2.4, $E_S = \varinjlim_t E_R \otimes_R
S/(\bz)\brt$.
\qed
\enddemo

\proclaim{Corollary 3.4} Let $(R,\m)$ be an excellent reduced local
ring.  Then the tight closure of $0$ in  $E_R$ is independent 
of whether one computes it over $R$ or over $\hatR$.  
\endproclaim

\demo{Proof}  The closed fiber of $R \to \Rhat$ is a field.  The
residue fields are the same; so the extension of residue fields
is separable.  By Lemma 2.2, $R$ and $\hatR$ share a common
CS test element.  Thus Theorem 3.3 applies. \qed
\enddemo

This shows that for $(R,\m)$ excellent reduced, $\ttau(R) = \Ann_R 0^*_{E_R} =
\ttau(\Rhat) \cap R$.  Thus $\ttau(R)\Rhat \inc \ttau(\Rhat)$.  Below we
show that equality does in fact hold.

By $(R,\um)$ we mean a semi-local ring with maximal ideals
$\um = \m_1,\ldots, \m_t$.  By completion we mean completion
at the Jacobson radical $\m_1 \cap \cdots \cap \m_t$.  This
is then a product of complete local rings and $E_R 
= \oplus_{i=1}^t E_{R_{\m_i}}$.
Note that when $(R,\um) \to (S,\un)$ is flat 
with Gorenstein closed fibers and $0^*_{E_S}$ is
extended from $0^*_{E_R}$, then $\ttau(R) = \ttau(S) \cap R$.

\proclaim{Remark 3.5}  {\rm If $(S,\m)$ is  excellent then the map $S \to \Shat$
has geometrically regular fibers.  This does not mean that for
$Q \in \Spec(\Shat)$ and $P = Q \cap S$ the map of residue fields
of $S_P \to \Shat_Q$ is separable.  However, if $Q$ is minimal over
$P\Shat$, then $\Shat_Q/Q\Shat_Q$ is a separable extension of $S_P/PS_P$.}
\endproclaim

The next theorem extends a result by Lyubeznik and Smith.  
\comment 
By
\defn{biequidimensional} we mean a ring $S$ such that for all maximal
ideals $\n \inc S$, $S_\n$ is an equidimensional local ring of the
same dimension as $S$.
\endcomment  

\proclaim{Theorem 3.6} Let $(S,\un)$ be a 
 semi-local reduced excellent ring.
Then
\roster
\item $\ttau(S)$ localizes properly at any prime $Q \in \Spec(S)$,
i.e., $\ttau(S)_Q = \ttau(S_Q)$.
\item $\ttau(\hatS) = \ttau(S) \hatS$.
\item 
 Let $(R,\m)$ be a  semi-local reduced excellent ring
 such that $(R,\um)\to (S,\un)$ is a flat semi-local
map 
with Gorenstein $F$-injective closed fibers and
such that if $\n$ is a maximal ideal of $S$ lying over $\m$ then
$R/\m \to S/\n$ is separable.  
\comment
Assume also that for any prime $P$
associated to $\ttau(R)$,
 the extension $R_P/PR_P \to S_P/PS_P$ is reduced
and (generically)
separable. 
\endcomment
 Assume that $R$ and $S$ have a common CS 
 test element.
 Then $\ttau(S) = \ttau(R)S$.
\endroster
\endproclaim

\demo{Proof}  If the above statements do not hold, then take
a counterexample with $S$ having the smallest possible dimension,
and take it with  the least of (1), (2), or (3) possible.  Note that
the statements hold for fields; so $\dim S > 0$.  Also, in each case
we may immediately localize at a maximal ideal of $S$ (and its contraction
in case (3)).

Suppose that the counterexample is to (1).  Then by the induction
hypothesis, for any chain of primes $Q \subsetneq P \subsetneq \n_i$
we have $\ttau(S_Q) = \ttau(S_P)_Q$.  Consider $S \to \hatS$.
By \cite{16}, $\ttau(\hatS)$ localizes properly at any multiplicative
system.  Let $P \in \Spec(S)$ have dimension 1 (i.e., $\dim S/P = 1$).
Then $\hatS_P$ is semi-local with maximal ideals precisely those primes
of $\hatS$ lying over $P\hatS$  (which are necessarily minimal
over $P\hatS$).  The induction hypothesis on (3)
gives $\ttau(S_P)\hatS_P = \ttau(\hatS_P)$, and by \cite{16},
$\ttau(\hatS_P) = \ttau(\hatS)_P$.  Thus
$$
\ttau(S_P) = \ttau(S_P)\hatS_P \cap S_P = (\ttau(\hatS))_P \cap S_P
= (\ttau(\hatS)\cap S)_P = \ttau(S)_P.
$$ 
(The last
equality follows from Corollary 3.4.)  This contradicts $S$ being
a counterexample.

Say the counterexample is to (2).  
Then $\ttau(S)\hatS \inc \ttau(\hatS)$,
and we want to see that this is an equality.  The associated primes
of the extended ideal are primes of $\hatS$ minimal among
those lying over the associated
primes of $\ttau(S)$.  Say that $Q$ is such a prime (lying over $P$).
  Suppose first
 that $Q$ is not maximal.  In this case we have
$(\ttau(\hatS))_Q = \ttau(\hatS_Q) = \ttau(S_P)\hatS_Q 
= \ttau(S)_P \hatS_Q = (\ttau(S)\hatS)_Q$.  The first
equality follows since $\ttau$ localizes properly
in complete rings, the second equality uses the induction hypothesis
(3) (and  Remark 3.5)
in a map of rings of smaller dimension, and the
third equality uses that $S$ has property (1).

We have reduced to the case that $\ttau(\hatS)/\ttau(S)\hatS$
is $0$-dimensional. So we may apply the following more general claim:
if $A \inc E = E_S(L)$, $I = (0:_S A)$, $J = (0:_{\hatS} A)$ and 
$\ell(J/I) < \8$ then $J = I\hatS$.  To see this there is no loss of
generality in modding out by $I$.  Thus $A$ is an $S$-faithful submodule
of $E$.  If $J \ne 0$ then $\ell(J) < \8$ implies that $\hatS$, and therefore
$S$, are depth $0$.  Hence $J$ is extended from $R$, and $A$ is not
faithful.  Thus $J$ must be $0$.

Suppose now that we have a counterexample to (3). 
 By Theorem 3.3,  $0^*_{E_S}$ is extended
from $0^*_{E_R}$.  
\comment 
Thus $\ttau(R) = \ttau(S) \cap  R$.  The associated primes
of $\ttau(R)S$ are primes minimal among those lying over associated primes
of $\ttau(R)$.  Let $Q
\in \Spec(S)$, lying over $P$, be associated to $\ttau(R)S$.  If
$Q \ne \n$ then applying the induction hypotheses yields
$$
\ttau(S)_Q = \ttau(S_Q) = \ttau(R_P)S_Q = 
(\ttau(R)_P)S_Q = (\ttau(R)S)_Q
$$
 and we are done.
Otherwise we may assume that $Q = \n$ and $P = \m$.  By part (2), 
we may complete both
$R$ and $S$.  The hypothesis now guarantees that
the  extension of residue fields is
separable, so $0^*_{E_S}$ is extended.
\endcomment
The following  lemma completes the proof. \qed
\enddemo

\proclaim{Lemma 3.7}  Let $(R,\m) \to (S,\n)$ be a flat
map of reduced  local rings such that $S/\m S$ is
Gorenstein and $0^*_{E_S}$ is extended.  Assume also that
$\ttau(\hatR) = \ttau(R)\hatR$ and $\ttau(\hatS) = \ttau(S)\hatS$. Then
$\ttau(S) = \ttau(R)S$.  If $0^{*fg}_{E_S}$ is extended, 
$\tau(\hatR) = \tau(R)\hatR$, and $\tau(\hatS) = \tau(S)\hatS$,
then
$\tau(S) = \tau(R)S$.
\endproclaim
\demo{Proof}  Let $\bz = \vec z t$ be elements of $S$ which
form an s.o.p.~for $S/\m S$.
Then we have 
$$
\split
\Hom_S(S/\ttau(R)S, E_S)  &= \Hom_S(R/\ttau(R) \otimes _R S,
\varinjlim_{t} E_R \otimes_R S/(\bz)\brt)  \\
& = \varinjlim_t \Hom_S(R/\ttau(R) \otimes S/(\bz)\brt, E_R \otimes
S/(\bz)\brt)  \\
& = \varinjlim_t \Hom_R(R/\ttau(R), E_R) \otimes S/(\bz)\brt \\
& = \varinjlim_t 0^*_{E_R} \otimes S/(\bz)\brt
= 0^*_{E_S}.
\endsplit
$$
The last equality follows from the assumption that
$0^*_{E_S}$ is extended.  The  second to last equality follows because
$\ttau(\hatR) = \ttau(R)\hatR$.  The other equalities follow from the flatness of
$S/(\bz)\brt S$ over $R$ and standard facts about direct limits.

Also, since $\ttau(\hatS) = \ttau(S)\hatS$, 
$0^*_{E_S} = \Hom_S (S/\ttau(S), E_S)$.
But now by applying the exact functor
$\Hom_S(-,E_S)$ to $0 \to C \to S/\ttau(R)S
\to S/\ttau(S) \to 0$ we obtain that $\Hom_S(C, E_S) = 0$.  Hence
$C = 0$. 

The last statement has the same proof.
\qed
\enddemo

Now that we have finished the proof that $\ttau$ extends properly, we
may omit that hypothesis from Lemma 3.7.

\proclaim{Corollary 3.8} Let $(R,\m) \to (S,\n)$ be a flat
map of reduced  local rings such that $S/\m S$ is
Gorenstein and $0^*_{E_S}$ is extended.   Then
$\ttau(S) = \ttau(R)S$.
\endproclaim
\demo{Proof}  By Theorem 3.6(2), the additional condition needed in
Lemma 3.7 is satisfied. \qed
\enddemo

In the situation that $(R,\m) \to (S,\n)$ is flat, we expect a common
test element when the generic fibers are Gorenstein and $F$-rational.
Below, we give a weaker condition which allows us to find a CS test
element for $S$ which is almost as useful as a common test element.

\proclaim {Theorem 3.9}  Let $(R,\m) \to (S,\n)$ be a  map
of excellent reduced rings.    Assume that
$\m S$ is reduced in $S$ and set $U \inc S$ to be
the complement of the minimal primes of $\m S$.  Suppose that  $R \to U\1 S$
is flat with $U\1 S/\m U\1 S$ separable over
$R/\m$. Then  there is an 
element $d \in U$  such that  
$d\ttau(R) \inc \ttau(S)$.  In particular $S$ has CS test elements of
the form $cd$, where $c$ is a CS test element of $R$ and $d \in U$.
\endproclaim

\demo{Proof}     
  The map
$R \to U\1 S$ is flat with Gorenstein
closed fibers of dimension $0$, and  
thus $E_{U\1 S} = E_R \otimes_R  U\1 S$ by Theorem 2.4.

Let $Q$ be a minimal prime of $\m S$ and set $T = \widehat{S_Q}$.
 The separability hypothesis  and 
Theorem 2.3 show that $\hatR$ and $T$ have a common $CS$ test element,
and thus by Theorem 3.3, $0^*_T$ is extended from $0^*_{\hatR}$.
Hence $0^*_{S_Q} = 0^*_T$ is extended from $0^*_{\hatR} = 0^*_R$.
By Theorem 3.6 and Corollary 3.8, $(\ttau(S))_Q = \ttau(S_Q) = \ttau(R)S_Q$
for each associated prime of $\m S$.  Thus there is a $d \in U$ such that
$d \ttau(R) S \inc \ttau(S)$.
\qed
\enddemo

We now state another version of Lemma 3.1.

\proclaim{Lemma 3.10}  Let $(R,\m) \to (S, \n)$ be a flat map of reduced local
rings with  $F$-rational closed fiber and also satisfying the conditions of 
Theorem 3.9. 
Suppose $\bz = \vec z d$ is an s.o.p. in $S/\m S$ and $b$ maps to
a socle element in $S/(\m, \bz)S$.  If $M$
 is an Artinian module and $u$ is not in $0^*_M$,
then $bu$ is not in $0^*_{M \otimes S/(\bz)}$.
\endproclaim

\demo{Proof}   Let $M' = M \otimes S/(\bz)$.  
Choose $c$ and $d$ as in Theorem 3.9.  If the lemma is not true,
then for all $q$,  $cd(bu)^q =0$ in $F^e_S(M')$.
Thus $db^q \in 0 :_S cu^q = (0 :_R c u^q)S + (\bz)\brq S$,
by flatness.  Since $u \notin 0^*_M$, $(0 :_R c u^q) \inc \m$
for infinitely many $q$.  Thus $db^q \in \left((\bz)\brq (S/\m S)\right)^*$, 
contradicting
the $F$-rationality of $S/\m S$. \qed
\enddemo

\head \S 4. Extension of tight closure by flat maps \endhead

We are now in a position to consider the extension of tight closure under flat
maps of local rings.  In Theorem 4.1 we consider the case where the closed
fiber is Gorenstein and either $F$-injective or $F$-rational.  In
Theorem 4.2 we weaken the hypothesis and assume only that the closed
fiber is CM.  As a trade-off, however, we must then assume that the closed
fiber is either geometrically $F$-injective or geometrically $F$-rational.
As an application, we extend the result (independently obtained by
Enescu \cite{7} and Hashimoto \cite{8}) that the extension of an $F$-rational
ring by a flat map with geometrically $F$-injective closed fiber is $F$-rational.

\proclaim{Theorem 4.1}  Let $(R,\m,K) \to (S,\n,L)$ be a flat
map of excellent rings
with Gorenstein closed fiber.  Assume that $K \to L$ is
a separable field extension.  Let $\bz = \vec z d \in S$ be elements
in $S$ which form an s.o.p.~in $S/\m S$.  Suppose that
\roster 
\item  $c$ is a common CS test element for $R$ and $S$,
and $S/\m S$ is $F$-injective, or
\item  
$S/\m S$ is $F$-rational and $K \to S/\m S$ is separable.

\endroster

Let $M$ be an $R$-module with DCC, and set $M' = M \otimes_R \varinjlim_t
S/(\bz)\brt S$.  Then $0^*_{M'}$ (respectively $0^{*fg}_{M'}$)
 is extended from $0^*_M$ (respectively $0^{*fg}_M$).
In particular, if $I$ is a tightly closed $\m$-primary ideal of $R$, then
$(I + (\bz)\brt)S$ is a tightly closed ideal of $S$ for all $t$.
\endproclaim

\demo{Proof}  
We may start by modding out by $0^*_M$ (respectively $0^{*fg}_M$)
 to assume that $0$ is
(finitistically) tightly closed in $M$.  Thus we wish to show that
$0$ is also (finitistically) tightly closed in $M'$.  
Note that since $S/\m S$ is Gorenstein, $\bz$ forms a regular sequence on $M \otimes_R S$. Hence,  the maps $M \otimes_R S/(\bz)\brt S \to
M \otimes_R S/(\bz)^{[t+1]}S$ are injective. So any socle element in $M'$ is
already present in $M \otimes S/(\bz)S$.

The initial setup here is as in Theorem 3.3.
Let $b \in S$ map to the socle of $S/(\m + (\bz))S$.  Given any 
choice of $\vec u n$ generating the socle of $M$, the set
$\{u_1 b,\ldots, u_n b\}$ generates the socle of $M'$
(the closed fiber is Gorenstein).
If $0$ is not tightly closed in $M'$, then there is a socle element
$v = \sum_{i=1}^n s_i u_i b \in 0^*_{M'}$, and the $s_i$ can
be considered to be elements of $L$.  Choose $\vec u n$ so that such 
a sum can be written with the fewest possible $u_i's$.  Let $h$ be
this number and reorder so that $v = \sum_{i=1}^h s_i u_i b$.
Suppose that the elements $\vec s h$ (thought of as elements of $L$)
are linearly dependent over $K$.  Then without loss of generality,
$s_h = r_1 s_1 + \cdots + r_{h-1} s_{h-1}$ and hence 
$v =  \sum_{i=1}^{h-1} s_i(u_i + r_i u_h)b$,
contradicting our choice of $h$.  Also, $h \ne 1$ by Lemma 3.1.

In both cases we can choose $c$ and $d$ as in Theorem 3.9
(where in case (1) we may assume that $d=1$).

Since $v \in 0^*_{M'}$ we have that $cdv^q = 0$ in $F^e_S(M')$. Using the definition of the direct limit and 
that $\bz$ is a regular sequence on $M \otimes_R S$, we get 
$cdv^q = 0$ in $F^e_S(S/(\bz)S \otimes_R M)$.

We have $d b^q
\in \left(c(\vec u {h-1})\brq :_S c u_h^q\right)$
where the elements defining the colon ideal belong to $F^e_S(S/(\bz)S \otimes_R M)$.
As we noted earlier, $R \to S/(\bz\brq)S$ is flat, and so

$d b^q \inc ((c(\vec u {h-1})\brq):_R cu_h^q)S + (\bz)\brq S$ for all $q$, where the elements defining the colon ideal
belong to $F^e_R(M)$.

If this  colon ideal is contained in $\m$  for $q \gg 0$, then we
get  a contradiction (see the cases listed below).
So we may assume that for large enough $q_1$, $c u_h^{q_1} 
= r_{h,1} cu_1^{q_1} + \cdots + r_{h,h-1} cu_{h-1}^{q_1}$ in $F^{e_1}_R(M)$.
Then $cd v^{q_1} = \sum_{i=1}^{h-1} cd (s_i^{q_1} + r_{h,i} s_h^{q_1})
u_i^{q_1} b^{q_1}$.  Each $s_i^{q_1} +
r_{h,i} s_h^{q_1}$ must be a unit, since $L$ is separable over $K$.
Thus we have $db^q \in  c(\vec u {h-2})\brq  :_S cu_{h-1}^q 
\inc ((c(\vec u {h-2})\brq):_Rc u_{h-1}^q)S + (\bz)\brq S$ for all $q \gg
 0$.  We obtain a contradiction (see the cases listed
below), unless for all $q_2 \gg 0$
 $cu_{h-1}^{q_2} 
= cr_{h-1,1} u_1^{q_2} + \cdots + cr_{h-1,h-2} u_{h-2}^{q_2} 
$ in $F^{e_2}_R(M)$.  Then
$cd v^{q_2} = \sum_{i=1}^{h-2} (s_i^{q_2} + r_{h,i}^{q_2/q_1} s_h^{q_2}
+ r_{h-1,i} s_{h-1}^{q_2})
u_i^{q_2} b^{q_2}$.  Again, each coefficient
must be a unit by separability.  Continuing in this manner we get
a unit $\al \in S$ and $q$ such that in the respective cases:
\roster
\item  $v^q = \al u_1^q b^q = 0$ in $F^e_S(M')$, implying that
$b^q \in (\bz)\brq (S/\m S)$, contradicting the $F$-injectivity of $S/\m S$ (look, for example, at
the proof of Lemma 3.1);
\item  $v^q = \al u_1^q b^q \in 0^*_{F^e_S(M')}$, contradicting
Lemma 3.10.
\endroster

The case of the finitistic tight closure of 0 is similar and we will leave it
as an exercise to the reader.
\qed
\enddemo

Now we would like to treat the case when the closed fiber is Cohen-Macaulay,
and not necessarily Gorenstein. 

\proclaim{Theorem 4.2}  Let $(R,\m,K) \to (S,\n,L)$ be a flat map
of excellent rings with CM closed fiber.  Let $I \inc R$
be $\m$-primary and tightly closed.  Let $\bz = \vec z d$ be 
elements whose image in $S/\m S$ is an s.o.p.  Suppose that either
\roster
\item the element $c \in \Ro$ is a common test element for $R$ 
and $S$, and that $S/\m S$ is geometrically 
 $F$-injective over $K$, or
\item $S/\m S$ is geometrically $F$-rational over $K$.
\endroster
Then $(I,\bz)S$ is tightly closed.
\endproclaim
\demo{Proof}  
Let $\{\vec b m\}$ map to generators of the socle of $S/(\m,\bz)S$.
If $(I,\bz)S$ is not tightly closed, then for any choice of generators
of the socle of $R/I$, say $\{\vec u n\}$, there exists a socle
element $w = \sum_{i,j} \la_{ij} u_j b_i$ in $((I,\bz)S)^*$.  Choose
$\{\vec u n\}$ and $w$ to minimize the number of nonzero $\la_{ij}$.
We first observe that $\{\la_{ij}b_i\}_{i,j}$ is part of a minimal 
basis for  $S/IS$ over $R/I$.  Otherwise, by Remark 3.2
there exist $k_{ij} \in K$ such that $\sum_{i,j} k_{ij} \la_{ij} b_i
= 0$ in $S/(\m,\bz)S$.  This would allow us to choose $\{u'_j\}$ to 
get a smaller number of nonzero $\la_{ij}$.  

We claim that in fact $\{\la_{ij}^qb_i^q\}_{i,j}$ is part of a
free $R/I\brq$-basis for $S/(I,\bz)\brq$ for every $q$.  
If not, for some $q$ there exist $k_{qij}
=k_{ij} \in K$
such that $\sum_{i,j} k_{ij} \la_{ij}^q b_i^q = 0$ in $S/(\m, (\bz)\brq)S$.
Let $K' = K[k_{ij}^{1/q}]$ and then the
element $(\sum_{i,j} k_{ij}^{1/q} \la_{ij} b_i)^q = 0$ in $S/(\m,(\bz)\brq)S
\otimes K'$.
But $S/ \m S\otimes K'$ is $F$-injective; so 
$\sum_{i,j} k_{ij}^{1/q} \la_{ij} b_i
=0$ in $S/(\m, \bz)S \otimes K'$, a contradiction.

In case (1), $c w^q = \sum_{i,j}(cu_j^q)\la_{ij}^q b_i^q \in (I,\bz)\brq$
then implies that $cu_j^q \in I\brq R$, or $u_j \in I^* = I$, a contradiction.

In case (2), choose $c \in \Ro$ and $d \in S^\circ$ as in Theorem
3.9.  We claim that  the set 
$\{d\la_{ij}^q b_i^q\}_{i,j}$ is part of a free $R/I\brq$-basis for
$S/(I,\bz)\brq$ for infinitely many $q$.
Otherwise, as in the analysis above we
 get  $d (\sum_{i,j} k_{qij}^{1/q} \la_{ij}
b_i)^q
\in (\bz)\brq S/\m S \otimes K^\infty$.  Each
$w_q = \sum_{i,j} k_{qij}^{1/q} \la_{ij} b_i$ is a socle element
in $S/(\m, \bd z) \otimes K^\infty$, and the socle is finite dimensional
as a $K^\infty$ vector space.
Let $V_q$ be the $K^\8$-subspace of the socle generated by
$\{w_{q'}\}_{q' \ge q}$.  Then by $DCC$ there exists $q_0$ such
that $V_q = V_{q_0}$ for all $q \ge q_0$.  Let $K' = K[k^{1/q_0}_{q_0 i j}]$
(this is a module-finite extension of $K$).  Then
$d w_{q_0}^q \in d(w_{q'}\mid q'\ge q_0)\brq  \cap S/\m S \otimes
K' \inc (\bd z)\brq S/\m S 
\otimes K'$.
  But this shows that a socle element is in the tight
closure of $S/(\m,\bd z) \otimes K'$,
 contradicting the 
$F$-rationality of $S/\m S \otimes K'$.
\qed
\enddemo

We can now improve Enescu's and Hashimoto's theorem on extension
of $F$-rational rings by flat maps with  sufficiently nice generic
and closed fibers.

\proclaim{Theorem 4.3}  Let $(R,\m, K) \to (S,\n, L)$ be a flat map
of excellent rings.  Assume that  $R$  is $F$-rational.  Suppose
that either 
\roster
\item $c \in \Ro$ is a common parameter test element for $R$ and $S$
and $S/\m S$ is geometrically $F$-injective, or
\item $S/\m S$ is geometrically $F$-rational.  
\endroster
Then $S$ is $F$-rational.
\endproclaim

\demo{Proof}
  If $R$
and $S$ have a common parameter test element, then case (1)
of Theorem 4.2 applies
with $I \inc R$ a parameter ideal (this case is known already by
the results of Enescu \cite{7} and of Hashimoto \cite{8}, under some finiteness conditions).  
Otherwise case (2) applies.
\qed
\enddemo

The next proposition generalizes \cite{9}, Theorem 7.36, where the closed
fiber is assumed to be regular.

\proclaim{Proposition 4.4} Let $(R,\m,K) \to (S,\n,L)$ be a flat map
of reduced excellent rings with Gorenstein closed fiber such that $R$ is
complete and either 
\roster \item  $S/\m S$ is $F$-injective, $K \subset L$ is separable
and $R$ and $S$ have a common test element; 
\item  $S/\m S$ is F-rational and $K
\to S/\m S$ is separable; 
\item the element $c \in \Ro$ is a common test element for $R$ 
and $S$, and  $S/\m S$ is geometrically 
 $F$-injective, or
\item  $S/\m S$ is geometrically $F$-rational over $K$.
\endroster

Then $\tau(R)S = \tau(S)$. \endproclaim

\demo{Proof} Both $R$ and $S$ are approximately Gorenstein. If ${I_{t}}$ is a
sequence of irreducible $\m$-primary ideals cofinal with  powers of $\m$, then
$\tau(R)= \cap I_{t}:I_{t}^*$.  Let $\bz = \vec z d \in S$ be an s.o.p.~for
$S/\m S$. Now, $J_{t}= (I_{t},(\bz)\brt)S$ is a sequence of irreducible ideals
in $S$ (because $S/\m S$ is Gorenstein). By Theorems 4.1 and 4.2, $J_{t}^*=
(I_{t}^*, (\bz)\brt)S$. Because of the flatness of the original map, we get that
$J_{t}:J_{t}^* = (I_t, \bz\brt)S:_S (I_t^*,\bz\brt)S
= (I_t^*:_R I_t)S + (\bz)\brt S$. Hence, 
$$ \multline \tau(S) = \cap_t (J_t:_S J_t^*) = \cap_t (I_t:_R I^*_t)S + (\bz)\brt
\\ = \cap_{t'}\left( \cap_t (I_t:_R I^*_t)S \right ) + (\bz)^{[t']})
= \cap_{t'} (\tau(R)S + (\bz)^{[t']}) = \tau(R)S
\endmultline$$ (in the next to last equality we use that $R$ is complete
and hence $S$ is $\cap$-flat over $R$).  \qed \enddemo

\proclaim{Remark 4.5}{\rm
Recently, Bravo and Smith have studied the behavior of test ideals under smooth and \'etale base change using different methods than ours (\cite {5}). Their results 
are stated under the assumption that the test ideal localizes or that the finitistic tight closure of zero equals the tight closure of zero in the injective hull of a reduced local excellent ring. These are hard open questions that are likely to be true under very mild conditions. Our results give some more  general conditions under which the test ideal extends naturally via flat base change. The assumption that $R$ is complete is not very restrictive. Comparing to the results of Bravo and Smith, if the test ideal localizes, it also extends well under completion. A proof of this can be given following closely the proof of Theorem
3.6 (one needs to use Theorem 4.1 instead of Theorem 3.3 whenever it is needed). Also, if the finitistic tight closure of 0 equals the tight closure of 0 in the injective hull, then the test ideal equals the CS test ideal which behaves well under completion (cf. Theorem 3.6).}
\endproclaim

\head \S 5.  Strong test ideals and applications
\endhead

An ideal $T \inc R$ is said to be a \defn{strong test ideal} for $R$
if $TI^* = TI$ for all ideals $I \inc R$ (see \cite{13}).  Vraciu
has shown that in a complete reduced ring $(R,\m,k)$, if $T$ is the annihilator
of any finitistically tightly closed submodule of $E_R(k)$, then
$T$ is a strong test ideal \cite{21}, Theorem 3.2.  As a corollary we
now obtain
\proclaim{Theorem 5.1}  Let $(R,\m,K)$ be an excellent reduced ring.
Then $\ttau(R)$ is a strong test ideal for $R$.
\endproclaim

\demo{Proof} This follows immediately from the faithful flatness of
$R \to \Rhat$ and from the fact that $\ttau(\Rhat) = \ttau(R)\Rhat$.
\qed
\enddemo

We would like to show that $\ttau(R)$ is a strong test ideal in the
CS sense,  i.e., if $M = \varinjlim_t R/I_t$ is Artinian, and
$w = \{w_t\} \in 0^*_M$, then for some $t$, $\ttau(R) w_t \in \ttau(R)I_t$.
We have not been able to prove this; however, we can prove a 
slightly weaker result which will suffice for our main application.

\proclaim{Theorem 5.2}  Let $(R,\m, K)$ be a reduced excellent ring.
Suppose that $M = \varinjlim_t R/I_t$, where each $I_t$ is $\m$-primary.
Let $w = \{w_t\} \in 0^*_M$.
Then for all $N >0$ there exists $t$ such that $\ttau(R) w_t \in
(\ttau(R)+\m^N)I_t$.
\endproclaim

\demo{Proof}  We give a proof here that is modeled after Vraciu's
proof of \cite{21}, Theorem 3.2.  As in the proof of Theorem 5.1,
it suffices to show that the statement holds when $R$ is complete.
Let $\ttau = \ttau(R)$. 

Fix $N >0$.  Let $M_1 = \varinjlim_t  R/(\ttau + \m^N)I_t$.  This
makes sense because for any ideal $A$ and any map $R/I \to R/J$
there is a clear induced map $R/AI \to R/AJ$.  The module $M_1$ is
supported only at $\m$.
Let $w' = \{w'_t\} \in M_1$ map to $w$ in $M$.  

We observe that $0:_{M_1} \ttau = \varinjlim_t (((\ttau +
\m^N)I_t):\ttau)/(\ttau + \m^N)I_t$ and clearly
$I_t \inc ((\ttau +\m^N)I_t):\ttau)$ and so 
$w' \in (0:_{M_1} \ttau)^*_{M_1}$.

Since $M_1$ is supported only at $\m$, we get an embedding $M_1 \to G$
where $G$ is
a (possibly infinite) direct sum of $E_R(k)$'s.  Then
$0:_{M_1} \ttau = M_1 \cap (0:_G \ttau)$.  Since $(0:_G \ttau)$ is
tightly closed, so is $0:_{M_1} \ttau$.  By the remark following
Proposition 2.1, $\ttau$ kills tight closure in all the modules
involved, and hence there exists $t$ such that
$\ttau w_t \in (\ttau +\m^N)I_t$.
\qed
\enddemo

The next theorem shows that tight closure in the injective hull may be
shown to extend without the need for a common CS test element as
long as the extension of residue fields is separable.  Because we
will make use of a result from \cite{3} (see Theorem 5.4 below)
which is proved for rings whose completions are domains, we will restrict
our statement to that case.  We believe, however, that Theorem 5.4 
is still true for reduced equidimensional rings, so that Theorem 5.3
should also hold under these more general assumptions.

\proclaim{Theorem 5.3} Let $(R,\m, K) \to (S,\n, L)$ be a flat map of
 excellent domains both of whose 
completions are domains.  Assume that $S/\m S$ is Gorenstein and $F$-rational
and $K \to L$ is separable.  Then $0^*_{E_S}$ is extended from $0^*_{E_R}$,
and $\ttau(S) = \ttau(R)S$.

Also, $0^{*fg}_{E_S}$ is extended from $0^{*fg}_{E_R}$ and hence $R$ and
$S$ have a common test element (we cannot make the assertion that
$\tau(S) = \tau(R)S$ from this).
\endproclaim

\demo{Proof}
Choose $\bz = \vec z d$ an s.o.p. for $S/\m S$ and let $b$
be the socle element in $S/(\m + \bz)$.  Let $\vec u h$ generate
the socle modulo $0^*_{E_R}$ where we write $E_R(k)$ as 
$\varinjlim_t R/I_t$.

Suppose now that $0^*_{E_S}$ is not extended.  Then there is an
element $ x = (\sum s_j u_j)b$ in $0^*_{E_S}$, i.e., there is
a $d \in S^\circ$ such that for all $q$ there exists
$t$ such that  $d(\sum s_j^q u_{j,t}^q)b^q \in (I_t,\bz)\brq S$.

We can choose $x$ and $\vec u h$ so that the number of nonzero terms in the expression
of $x$ is minimal. As in the proof of Theorem 4.1, it follows that the coefficients $s_{j}$
are linearly independent over $K$ when regarded as elements of $L$.

We can also arrange that for some $w$, $u_{w+1},\ldots, u_h \in
(u_1,\ldots, u_w)_{E_R}^*$, but $u_w \notin (u_1,\ldots, u_{w-1})_{E_R}^*$.
Also, suppose that $\ttau(R) = (\vec c l)$ minimally.

As in the proof of Theorem 4.1, we can easily show that $bu_{w} \notin 
(u_1,\ldots, u_{w-1})_{E_S}^*$ in $E_{S}$.

Fix $N > 0$. According to Theorem 5.2, for every $q$ we can find $t$ such that 
$$\ttau(R) (u_{t,w+1},\ldots, u_{t,h})\brq \inc
(\ttau(R)+\m^N) (u_{t,1},\ldots, u_{t,w}, I_t)\brq.$$

Multiply now by $c_1$ and use Theorem 5.2 to get for each $q$, an integer $t$
such that (omitting the $t$ indices on the left) 
$$ 
d \left( \sum_{j = 1}^w
(c_1 s_j^q + \sum_{v=w+1}^h s_{v}^q r_{j,v,q})u_j^q \right) b^q \in (I_t,
\bz)\brq S.
$$ 
Here the $r_{j,v,q}$'s are in $\ttau(R) +\m^N$. We then have 
$$
\split d (c_1 s_w^q + \sum_{v=w+1}^h s_v^q r_{w,v,q})b^q &\in ((u_1,\ldots,
u_{w-1},I_t, \bz)\brq :_S u_w^q ) \\
& = \left( (\vec u {w-1}, I_t)\brq:_R u_w^q\right)S + (\bz)\brq S
\endsplit  
$$  
and, since $u_w \notin
(u_1,\ldots, u_{w-1})_{E_R}^*$, by Theorem 5.4 below there is a 
$q_0$ such that 
$$
(u_1,\ldots, u_{w-1}, I_t)\brq:_R  u_w^q  \inc
\m^{[q/q_0]}.
$$  

Note that $b^{q_0} \notin ((\m + (\bz)^{[q_0]})S)^*$ (or else
$b \in ((\bz)S/\m S)^*$) so 
$$
\left( c_1(s_w^q + \cdots + r'_{h,q} s_h^q) + c_q\right) \in (\m
+(\bz)^{[q_0]})S:_S d(b^{q_0})^{q/q_0} \inc \n^{[q/q_1]}
$$
for some $q_1 \ge q_0$.
The element $c_q$ is in $(c_2,\ldots, c_l)+\m^N$.  
By separability of $L$ over $K$, the coefficient of $c_1$ must be a
unit for all $q$.
Thus applying the Krull intersection theorem we see that $c_1$ is in
$$
\left(\bigcap_q ((c_2,\ldots, c_h) +\m^N + \n^{[q/q_1]})\right) \cap R  
\inc ((c_2,\ldots, c_h)+\m^N)S \cap R = ((c_2,\ldots, c_h)+\m^N)R.
$$ 
This argument is valid for all $N >0$. So another application of
the Krull intersection theorem shows that $c_1 \in (c_2,\ldots, c_l)$,
a contradiction.  By Corollary 3.8, $\ttau(S) = \ttau(R)S$.  

The proof for the finitistic tight closure is exactly the same,
although we may use Vraciu's original result, since we may immediately
fix a value of $t$.  \qed
\enddemo

The next Theorem appears as \cite{3}, Proposition 2.4.  We have
altered the statement slightly to include the element $d \in \Ro$,
but the proof is essentially the same.

\proclaim{Theorem 5.4} {\rm }  Let $(R,\m)$ be an
excellent local domain such that its completion is a domain.  Let $M =
\varinjlim_t R/I_t$ be a direct limit system.  Fix $u \notin 0^*_M$
and $d \in \Ro$.  Then
there exists $q_0$ such that   $J_q = \cup_t (I_t\brq:d u_t^q) \inc \m^{[q/q_0]}$
for all $q \gg 0$ (where $\{u_t\}$ represents $u \in M$ and $u_t \mapsto
u_{t+1}$).   In particular, if $I \inc R$ we may take $M = R/I$ where the limit
system consists of equalities. Then $u \notin I^*$ implies that $(I\brq:du^q)
\inc \m^{[q/q_0]}$. \endproclaim

\proclaim{Remark 5.5}  {\rm As mentioned above, Theorem 5.3 should be true
whenever $R$ and $S$ are reduced and equidimensional (and excellent).
The proof given above does give the result when the completion of
$R$ is a domain, the closed fiber of the map is regular,
$S$ is reduced (not necessarily equidimensional) and $K \inc L$ is separable.  This is
because we may pick $\bz$ to be a regular system of parameters
for $S/\m S$ and take $b = 1$. When $b=1$, using that there is $q_0$ such that
$(u_1,\ldots, u_{w-1})\brq:_R du_w^q + I_{t}^{[q]} \inc \m^{[q/q_{0}]}$, one has that  $c_1(s_w^q + \cdots + r'_{h,q} s_h^q) + c_q \in
(\m + (\bz)^{[q_0]})^{[q/q_0]}$, which is clearly in
larger and larger powers of $\n$ as $q$ increases.  The rest
of the argument then follows.}
\endproclaim

We can now give an improved version of Theorem 3.9.

\proclaim{Theorem 5.6} Let $(R,\m) \to (S,\n)$ be a map of
excellent reduced rings.  Assume that $\Rhat$ is
a domain.  Let $Q \in \Spec(S)$ be such that $\m S \inc Q$,
$R \to S_Q$ is flat, $\m S_Q$ is reduced, and $R/\m \to S_Q/QS_Q$ is separable.  Assume that either
\roster
\item
$\widehat {S_Q}$ is a domain and $\widehat {S_Q}/ \m \widehat {S_Q}$
is Gorenstein and $F$-rational, or
\item
 $S_Q/\m S_Q$ is regular.
\endroster
Then there is an element $d \in S -Q$ such that $d \ttau(R) \inc \ttau(S)$.
\endproclaim

\demo{Proof}  Given Theorem 5.3 (or Remark 5.5), the proof is
identical to that of Theorem 3.9.
\qed
\enddemo

We will end with several more consequences of the ideas used in the proof of
Theorem 3.3. The first one addresses the notion of Hilbert-Kunz multiplicities.
For an $\m$-primary ideal $I$ in a local ring $(R, \m, K)$ of Krull dimension $d$, the
Hilbert-Kunz multiplicity of $I$ is defined as  $e_{HK}(I,R) = \lim_{q
\rightarrow \infty} l_R(R/I\brq)/q^d$. If $I = \m$, we will just write
$e_{HK}(R)$.

\proclaim{Proposition 5.7} Let $(R,\m,K) \to (S,\n,L)$ be a flat local map with a
Cohen-Macaulay closed fiber. Then for every $\m$-primary ideal $I$ and
for every s.o.p. $\bz$ of
the closed fiber $S/\m S$,
$$e_{HK}(I,R) e_{HK}(S/\m S) = e_{HK}(IS+\bz S, S).$$ 
\endproclaim

\demo{Proof}
As before, 
$R/IR \to S/(IS+\bz S)$
is flat. Apply Remark 3.2 and deduce that in fact $S/(IS+\bz) S$ is free over $R/I$ of rank equal to the rank of $S/(\m S+\bz S)$ over $K$. We can apply this for $I^{[q]}$ and $z^{[q]}$. Denote by $d$ and $e$ the Krull dimension of $R$ and $S/\m S$, respectively. In fact, it follows that
$$l_S(S/(I\brq, \bz\brq)S) = l_{R}(R/I^{[q]})l_{S}(S/(\m S+\bz\brq)S )$$ (one can see this by applying \cite{6}, Exercise 1.2.25).
Hence, 
$$(1/q^{(d+e)})  l_{S}(S/(I^{[q]}+\bz\brq)S) = (1/q^d) l_{R}(R/I^{[q]}) (1/q^e)l_{S/\m S}(S/(\m +\bz\brq)S.$$ 

Taking the appropriate limits we get the stated identity.
\qed
\enddemo

It might be helpful to note the following 

\proclaim{Corollary 5.8 (Kunz)} In the situation stated above, $e_{HK}(S) \leq
e(S/\m S)e_{HK}(R)$. Here, $e(S/\m S)$ stands for the multiplicity of $S/ \m
S$. \endproclaim

\demo{Proof}
Take $I=\m$ and use the fact that $e_{HK}(S) \leq e_{HK}(\m S+zS)$ and that $e_{HK}(S/\m S) \leq e(S/\m S)$.
\qed
\enddemo

This Corollary has been obtained earlier by Kunz in 1976 with a different
proof.

Lastly we would like to show how this circle of ideas may be applied
to the theory of test exponents, a notion introduced recently by Hochster and
Huneke \cite{12}. The notion is intimately connected to the localization problem in tight
closure theory.

First, we would like to remind the reader of the definition of a test exponent.
Let $R$ be a reduced Noetherian ring of positive characteristic, $c$ a fixed
test element for $R$, and $I$ an ideal of $R$. We say that $Q$ is a test exponent
for $c,I$ if whenever $cu^q \in I^{[q]}$ and $q \geq Q$, then $u \in I^*$. The
minimal test exponent for $c,I$ is denoted by $\exp(c,I,R)$, if it exists.

\proclaim{Proposition 5.9} Let $(R, \m,K) \to (S,\n,L)$ be a flat local map with
Gorenstein F-injective closed fiber. Assume that $K \subset L$ is separable.
Let $c$ be a common completely stable test element of $R$ and $S$. If a test
exponent exists for $c^2, I, R$, then a test exponent exists for $c,  IS+\bz S,
S$. Moreover, $\exp(c^2,I,R)  \geq \exp(c, IS+\bz S, S)$. Here, $\bz$ stands for
a system of parameters of the closed fiber $S/\m S$. 
\endproclaim

\demo{Proof} Denote by $Q$ the test exponent of $c^2$ with respect to $I$ and
$R$. We will freely use the notation of Theorem 3.3 for the case $M =
R/I$. We want to show that if $cv^q \in I^{[q]}S+(\bz)^qS$ (let us denote this
relation by (\#)) for some $q \geq Q$, then $v \in (IS+\bz S)^{\ast}$. We know
that $(IS+\bz S)^{\ast} = (I^{\ast}S+\bz) S$. Fix $q \geq Q$ and denote by $J$
the ideal of elements $v$ that satisfy (\#). Clearly $(I^{\ast}S+\bz S) \subset
J$. We claim that these ideals are actually equal. If not, take $v \in J$ which
is in the socle of $S/(I^{\ast}S +\bz S)$. As in the proof of Theorem 3.3, we
write $v = \sum s_{i}u_{i}b$, where $u_{i}, i=1,...,n$ are generators for the
socle of $R/I^{\ast}$. As in Theorem 3.3, we get that $cu_{i}^Q \in (I^{\ast})^{[Q]}$.
But $c$ is a test element, so $c^2u_{i}^Q \in I^{[Q]}$. Therefore, $ u_{i} \in
I^{\ast}$. So, $v \in I^{\ast}S$, which is a contradiction, and we are done.
\qed  \enddemo

\enddocument